 \documentclass[a4paper,12pt]{article}
\usepackage{amssymb}
\usepackage{epsfig}
\usepackage{latexsym}
\usepackage{amscd}
\usepackage{amsmath}
\usepackage{amsthm}

\newcommand{\bC}{\ensuremath{\mathbb{C}}}

\newcommand{\bI}{\ensuremath{\mathbb{I}}}

\newcommand{\bR}{\ensuremath{\mathbb{R}}}

\newtheorem{theorem}{Theorem}[section] 
\newtheorem{lemma}[theorem]{Lemma}     
\newtheorem{corollary}[theorem]{Corollary}
\newtheorem{proposition}[theorem]{Proposition}
\newtheorem{remark}{Remark} [section]
\newtheorem{definition}{Definition}[section]

\begin{document}
\begin{center}
{{\Large{\bf 
Wave fronts and (almost) free divisors 
} }
}
\end{center}

\vspace{1pc}  \begin{center}{\large{ Susumu Tanab\'e}}\end{center}   %
\begin{center}
{{\it Dedicated to 60th birthday of David Mond}}
\end{center}
%
\noindent
\begin{center}
 \begin{minipage}[t]{10.2cm}
{\sc Abstract.} {\em  In this note we present a description of  wave front evolving from an algebraic hypersurface
 by means of a pull-back of the discriminantal loci of 
a tame polynomial via a polynomial mapping.  As an application we give examples of  wave fronts which define free/almost free divisors near the focal point.}
\end{minipage} \hfill
\end{center}
 \vspace{1pc}


\section{Introduction} 
\label{intro}

\noindent 
During the last decade we witnessed an intensive development of studies on wave fronts 
mainly from the differential geometric point of view.
We restrict ourselves to recall the works of authors like S.Izumiya, K.Saji, M.Umehara, K.Yamada.
After their definition, a map germ $f: (\bR^n,0) \rightarrow (\bR^{n+1},0)$ 
is called a wave front if there exists a unit vector field $\bf e$ along $f$
such that the $(f, {\bf e}): (\bR^n,0) \rightarrow (PT \bR^{n+1},0)$
is a Legendrian immersion.

In this article we propose an approach to the geometric studies of wave front (equidistant, parallel) surfaces based on complex analytic 
tools. The main difference between our method and the above mentioned differential geometric method consists in the way to represent a wave front: 
we try to describe it by means of its defining equation, while the latter relies on its parametrisation.
Here we present a  description of  wave front evolving from an algebraic hypersurface
 by means of the discriminantal loci of  a tame polynomial. 
In contrast to \cite{Ar76}, we do not identify the deformation parameter space of a polynomial and the space-time variables $(x,t).$
This situation makes us to consider a mapping $\iota$ $(1.7)$ from the physical space-time $(x,t)$ world to the deformation parameter space.
As we treat the deformation of the phase function $\Psi(x,t,z)$ $(1.4)$ (or square distant function $\Phi(x,t,z)$, Remark ~\ref{remark1})
in a global setting, we are obliged to  take into account its vanishing cycles at infinity.

\vspace{1.5pc}
\footnoterule

 \footnotesize{
AMS Subject Classification:  14B07, 14B05.
 \footnotesize{Partially supported by JSPS grant in aid (C) No.20540086}

\normalsize
To avoid difficulties
associated to the vanishing cycles at infinity with non-trivial monodromy (denoted by $E_c'$ in  \cite[Theorem 0.5]{DiS}), we have chosen a strategy
to treat only tame polynomial cases (see Lemma~\ref{lemma2}).

The majority of ever existing works reduce the Lagrangian or Legendrian singularities to local normal forms with the aid of diffeomorphisms.
In  \cite{Tan99} we have proposed an asymptotic analysis of the fundamental solution to hyperbolic Cauchy problem around the 
singular loci of the globally evolving wave front.
In this article we removed the quasihomogeneity condition imposed on the initial wave front in \cite{Tan99}.
Our Theorems   ~\ref{theorem2}, ~\ref{theorem3} generalise  \cite[Theorem 10]{Tan99}.

Another objective of this article is to show that certain wave front gives example of an almost free divisor.
The research on this kind of divisor has been initiated by J.Damon~\cite{Damon} and D.Mond~\cite{Mond}. The former gave the rank of "singular vanishing cycles" module while the latter uses differential forms to describe the vanishing cohomology in the singular Milnor fibre of an almost free divisor. 
 Their activities are motivated by a discovery \cite{DM} on the stabilization of a singular mapping: the discriminant of a stabilization 
 plays a role of "Milnor fibre" for the discriminant of the mapping.
 Despite remarkable results on topological invariants of almost free divisors, quite few non-trivial examples
are present in their works, especially those with the physical meanings are absent. Here we supply a class of examples that arise from geometric optics.
We shall notice, however, our basic idea to construct an almost free divisor  $\iota^{-1}(D_\varphi)$ as a pull back of a free divisor $D_\varphi$ (the definition of its freeness as a divisor coincides with the formulation oheorem ~\ref{theorem2}, 1) below) via a polynomial mapping
has already been underlined by J.Damon in \cite[p.219]{DM} and ~\cite{Damon} where he uses the terminology "nonlinear section" of a free divisor.

It is a great pleasure for us to dedicate our humble result to the 60th birthday of David Mond who founded a beautiful and significant theory on interrelations between free and almost free divisors. This relation would have never come into author's attention if
not twice stays at University of  Warwick realised by his invitation.

\section{Preliminaries on the wave fronts}

In this section we prepare fundamental notations and lemmata to 
develop our studies in further sections.
Let us denote by $Y:=\{(z,u) \in \bC^{n+1}; F(z)+u=0  \}$ the complexified initial wave front set
 defined by a polynomial  $F(z) \in \bR[z_1,\cdots, z_{n}] $, $z=(z_1,\cdots, z_{n})$. 
The real initial wave front set is given by
$Y \cap \bR^{n+1}$.

Let us consider the traveling of the ray starting from a point $(z,u) \in Y$ along unit vectors perpendicular to the hypersurface
tangent to $Y$ at $(z,u).$ It will reach at the point $(x_1, \cdots, x_{n+1})$
$$ x_j = \pm t \frac { 1}{| (d_z F(z),1)|} \frac{\partial F(z)}{\partial z_j} + z_j, 1 \leq j \leq n,$$
$$x_{n+1} = \pm t \frac {1}{| (d_z F(z),1)|} + u  \;\;{\rm with } \;\; (z,u) \in Y, \eqno(1.1)$$
 at the moment $t.$
Further on, we denote by $x'=(x_1, \cdots, x_n),$ $x=(x',x_{n+1}).$
 We see that $(x,t)$ and $(z,u)$ satisfying the relation $(1.1)$ are located on the zero loci of two phase functions
$$\psi_{\pm} (x,t,z,u) = \left( \langle x'-z, d_z F(z) \rangle  \ + (x_{n+1}-u) \right) \pm t| (d_z F(z), 1) |,
\eqno(1.2) $$ each of which corresponds to the backward $\psi_{+} (x,t,z,u)$
(resp. the forward $\psi_{-} (x,t,z,u)$ ) wave propagation.
To simplify the argument, we will not distinguish forward and backward wave
propagations in future. This leads us to introduce an unified phase function
$$\psi(x,t,z,u) := \psi_{+} (x,t,z,u) \cdot \psi_{-} (x,t,z,u)$$
$$=  \left( \langle x'-z, d_z F(z) \rangle  \ + (x_{n+1}+u) \right)^2 - t^2| (d_z F(z), 1) |^2, \eqno(1.3) $$

Let us denote by  $W_t$ the wave front at time $t$ with the initial wave front $Y$ i.e. $Y=W_0.$
Now we consider the following projection.
\begin{center}
$
\begin{array}{ccc}
\pi: \{(z,u) \in  Y: \psi(x,t,z,u)=0 \} & \rightarrow & \bC^{n+2}\\
  (x,t,z,u) & \mapsto & (x,t).
\end{array}  
$
\end{center}

\begin{lemma}
For  $x \in W_t$, the point $(x,t)$ belongs to the critical value set of the projection $\pi$ defined just above.
\label{lemma1}
\end{lemma}

We can understand this fact in several ways. Instead of purely geometrical interpretation, in our previous publication
\cite{Tan99}
we adopted investigation of  the singular loci of the integral of type, 
$$ 
I(x,t) = \int_{\gamma}  H(z,u)( \frac{1}{\psi_{+}(x,t,z,u)} + \frac{1}{\psi_{-}(x,t,z,u)})dz \wedge du
 $$
 for $\gamma \in H_{n}(Y)$ and $H(z,u) \in {\mathcal O}_{\bC^{n+1}}$. The above integral ramifies around
 its singular loci $W_t$ and by the general theory of the Gel'fand-Leray integrals (cf. ~\cite{Vas2}), 
$W_t$ is contained in the critical value set
 mentioned in the Lemma~\ref{lemma1}.

According to the Lemma ~\ref{lemma1},
The set $BW:=\cup_{t \in \bC} W_t \subset \bC^{n+1}$ (the real part of it is the big wave front after Arnol'd ~\cite[6.3]{Ar76}, ~\cite[22.1]{AVG})
can be interpreted as a subset of the discriminant
of the function (called the phase function)
$$ \Psi(x,t,z):=\left(\langle x'-z, d_{z} F(z) \rangle+ x_{n+1}+F(z) \right)^2-t^2 (| d_{z} F(z) |^2+1) \eqno(1.4)$$
for $x'=(x_1, \cdots, x_{n}).$  This is a set of $(x,t)$ for which the algebraic variety
$$X_{x,t}:=\{z \in \bC^{n}:  \Psi(x,t,z)=0 \}$$ 
has singular points.

\begin {remark}
Masaru Hasegawa and Toshizumi Fukui ~\cite{Hase} study the wave front $W_t$  as a discriminantal loci of the function,
$$ \Phi(x,t,z) =-\frac{1}{2}( |(x'-z,x_{n+1}+F(z)) |^2-t^2),$$
that measures the tangency of the sphere $\{ (z,z_{n+1}) \in \bR^{n+1}:|(z-x',z_{n+1}-x_{n+1}) |^2=t^2 \}$
with the hypersurface $Y \cap \bR^{n+1}.$ In some cases, this approach allows us to get less complicated 
expression of the defining equation of $BW$ in comparison with ours in Theorem ~\ref{theorem2}.

As Hasegawa points out, generally speaking, the inclusion of  $BW$ into the discriminant of $(1.4)$ is strict. The parabolic points of $F(z)$ produce
additionally so called "asymptotic normal surface." The discriminant, however, represents the big wave front in the neighbourhood of 
a focal point (see below).  Hence in our further studies on the local property around the focal point, this difference is negligeable.
\label{remark1}
\end{remark}

We assume that the variety $X_{x,t}$ has at most isolated singular points for a point $(x,t)$
of the space-time. Among those points, we choose a focal point $(x_0,t_0) \in \bC^{n+2}$ 
i.e. the point where the maximum of the sum of all local Milnor numbers is attained.
If we denote by $z^{(1)}, \cdots , z^{(k)}$ the singular points located on $X_{x_0,t_0}$
and Milnor numbers corresponding to these points by 
$\mu(z^{(i)}),$ $i=1,...,k$,   the following inequality holds for the focal point
$$ {\rm{sum \; of\; Milnor\; numbers \; of \; singular\; points\; on}}\; X_{x,t} \leq \sum_{i=1}^k \mu(z^{(i)}),$$
for every $(x,t) \in \bC^{n+2}.$  

Assume that the quotient ring
$$ \frac{\bC[z]}{( d_z \Psi(x_0,t_0,z) ) \bC[z]}  \eqno(1.5)$$
is a $\mu$ dimensional $\bC$ vector space  that admits a basis $\{e_1(z),$ $\cdots,$ $e_{\mu}(z)\}$
that contains a set of basis elements as follows,
$$e_1(z)=1, e_{j+1}(z)=(z_j-z_j^{(i)}), 1 \leq j \leq n, \eqno(1.6)$$ for a fixed
$ i \in [1,k].$ Here we remark that $\sum_{i=1}^k \mu(z^{(i)}) \leq \mu$. The denominator 
$( d_z \Psi(x_0,t_0,z) ) \bC[z]$ of the expression $(1.5)$ means the 
Jacobian ideal of the polynomial $\Psi(x_0,t_0,z).$

Now we decompose the  difference
$$\Psi(x,t,z)- \Psi(x_0,t_0,z) = \sum_{j=1}^m s_j(x,t) e_j(z)$$
by means a set of polynomials in $z$, $\{e_1(z), \cdots, e_{\mu}(z), e_{\mu+1}(z), \cdots, e_m(z) \}$
and a set of polynomials in $(x,t)$,

\begin{center}
$
\begin{array}{cccc}
\iota:\bC^{n+2} &\rightarrow &\bC^m\\
 (x,t) & \mapsto &\iota(x,t):=( s_1(x,t), \cdots, s_m(x,t))
\end{array} $
\end{center}
$$\eqno(1.7)$$
thus defined.In this way we introduce a set of polynomials $\{e_{\mu+1}(z),$ $\cdots,$ $e_m(z)$ $\}$
in addition to the basis of $(1.5).$
We consider a polynomial $\varphi(z,s)$ $\in$ $\bC[z,s]$ for $s=(s_1, \cdots,s_m)$
defined by 
$$\varphi(z,s)= \Psi(x_0,t_0,z) +\sum_{j=1}^m s_j e_j(z).  \eqno(1.8)$$

Locally this is a versal (but not miniversal) deformation of the holomorphic function germ
$\Psi(x_0,t_0,z)$ at $z=z^{(i)}.$

\section{Discriminant of a tame polynomial}

\begin {definition}
The polynomial $f(z) \in \bC[z]$
is called tame if there is a compact neigbourhood $K$ of the critical points of $f(z)$ such that 
$ \| d_z f(z) \| = \sqrt{(d_z f(z), \overline{d_z f(z)})} $ is away from 0   for all $z \not \in K.$
\label{def1}
\end{definition} 

In the sequel we use the notation $s'=(s_2,\cdots,s_m)$ and $s=(s_1,s').$

{Further on we impose the following conditions on $\varphi(z,s)$ introduced in $(1.8).$
Assume that there exists an open  set $ 0 \in V \subset \bC^{m-1}$ such that
$$ dim_{\bC} \frac{\bC[z]}{( d_z \varphi(z,s) ) \bC[z]}   < \infty, \eqno(2.1) $$
for every $s' \in V$ and $s_1 \in \bC$.
In addition to this, we assume that 
for every $s=(s_1,\cdots, s_{n+1}, 0,\cdots,0) \in \bC \times V,$ the equality 
$$ dim_{\bC}\frac{\bC[z]}{\left(d_z(\Psi(x_0,t_0,z)+\sum_{j=2}^{n+1} s_j e_j(z))\right) \bC[z]} =\mu, \eqno(2.2)$$
holds.
\begin {lemma}
Under the conditions $(1.5), (2.1),$ $(2.2)$ there exists a constructible subset ${\tilde U} \subset V$,
such that $\varphi(z,s)$ is a tame polynomial for every $s \in \bC \times {\tilde U}$ and
$$ dim_{\bC}\frac{\bC[z]}{( d_z \varphi(z,s) ) \bC[z]}  =\mu, $$
for every $s \in  \bC \times {\tilde U}$.
\label{lemma2}
\end{lemma} 
{\bf Proof} By \cite[Proposition 3.1]{Br88} $(2.2)$ yields the tameness of $\varphi(z,0).$
After Proposition 3.2 of the same article, the set of $s$
such that $\varphi(z,s)$ be tame is a constructible subset (i.e. locally closed set with respect to the Zariski topology) of the form $\bC \times W $ for $ W \subset V$.
According to \cite[Proposition 2.3]{Br88} the set 
$$T_n=\{ s \in \bC \times W: dim_{\bC} \frac{\bC[z]}{( d_z \varphi(z,s) ) \bC[z]} \leq n \},$$
is Zariski closed for every $n.$
We can take  $\bC \times {\tilde U} = T_{\mu} \setminus T_{\mu-1}.$
{\bf Q.E.D.}

{\bf Assumption I}

(i) By shrinking ${\tilde U}$ if necessary, we assume that a constructible set $U \subset {\tilde U}$ can be given 
 locally by holomorphic functions $(s_{\nu+1}, \cdots, s_m)$ on the coordinate space with variables 
$(s_2, \cdots, s_\nu),$ $\nu \geq \mu$.

(ii) The image of the mapping $\iota$ of a neighbourhood  of $(x_0,t_0)$ is contained in $\bC\times U.$ In other words,
$$ \iota (\bC^{n+2}, (x_0,t_0)) \subset (\bC\times U, \iota(x_0,t_0)).$$

For a fixed $\tilde s'=(\tilde s_2, \cdots, \tilde s_m)\in U$ and the constructible subset $U \subset V$ of 
the Assumption I,(i) we see that $\varphi(z, s_1, \tilde s')$ 
is a tame polynomial for all $s_1 \in \bC.$ For such  $\varphi(z, s_1, \tilde s'),$ 
we define the following modules,

$${\mathcal P}_\varphi(\tilde s') := \frac{\Omega_{\bC^n}^{n-1}}{d_z\varphi(z, s_1, \tilde s')
\wedge  \Omega_{\bC^n}^{n-2} +  d \Omega_{\bC^n}^{n-2}}, \eqno(2.3)$$

$${\mathcal B}_\varphi(\tilde s') := \frac{\Omega_{\bC^n}^{n}}{d_z\varphi(z, s_1, \tilde s')
\wedge  d\Omega_{\bC^n}^{n-2}}.\eqno(2.4)$$
the module ${\mathcal B}_\varphi(\tilde s')$ is called an algebraic Brieskorn lattice.
In considerig the holomorphic forms multiplied by $\varphi(z, s_1, \tilde s')$ be zero in $(2.3),$ $(2.4)$ we can treat two modules as $\bC[s_1]$ modules.

These modules contain the essential informations on the topology of the variety
$$ Z_{(s_1, \tilde s')} =\{z \in \bC^n : \varphi(z, s_1, \tilde s')=0\}. \eqno(2.5)$$ 
Let us denote by $D_\varphi \subset \bC \times U$ the discriminantal loci of the polynomial 
 $\varphi(z, s)$ i.e.
 $$ D_\varphi := \{s \in \bC \times U : \exists z \in Z_s, \; s.t. \; d_z \varphi(z, s)=\vec {0} \}. \eqno(2.6)$$

\begin{theorem} For a fixed $\tilde s'=(\tilde s_2, \cdots, \tilde s_m)\in U$,
both 
${\mathcal P}_\varphi(\tilde s')$ and ${\mathcal B}_\varphi(\tilde s')$
are free $\bC[s_1]$ modules of rank $\mu.$
\label{theorem1}
\end{theorem}
{\bf Proof} First we show the statement on ${\mathcal B}_\varphi(\tilde s')$.
After \cite[Theorem 0.5]{DiS} the algebraic Brieskorn lattice  ${\mathcal B}_\varphi(\tilde s')$ 
is isomorphic to a free $\bC[s_1]$ module of finite rank ( so called the Brieskorn-Deligne lattice). The topological triviality of the vanishing cycles at infinity for 
$\varphi(z, s_1, \tilde s')$  ensures this isomorphism.

On the other hand, for $(\tilde s_1, \tilde s') \in \bC \times U,$ the Corollary 0.2 of the same article tells us the following equality. 
$$  dim\;Coker (s_1-\tilde s_1|{\mathcal B}_\varphi(\tilde s') ) $$
$$= dim \;H_{n-1}(Z_{(\tilde s_1, \tilde s')})+ {\rm{sum \; of\; Milnor\; numbers \; of \; singular\; points\; on}}\;Z_{(\tilde s_1, \tilde s')}. $$
For $(\tilde s_1, \tilde s') \in \bC \times U \setminus  D_\varphi,$ the right hand side of the above equality equals 
$$ \sum_{s_1: Z_{( s_1, \tilde s')}\; \rm{singular}} {\rm{sum \; of\; Milnor\; numbers \; of \; singular\; points\; on}}\;Z_{( s_1, \tilde s')} $$
 by \cite[Theorem 1.1, 1.2]{Br88}.

Now we show that ${\mathcal B}_\varphi(\tilde s')$ is isomorphic to ${\mathcal P}_\varphi(\tilde s')$.

We show the bijectivity of the mapping $d: {\mathcal P}_\varphi(\tilde s') \rightarrow {\mathcal B}_\varphi(\tilde s').$
To see the injectivity, we remark that  the condition $d(\omega+ d\alpha+ \beta \wedge d\varphi(z, s_1, \tilde s')) = d \omega+ d\beta \wedge d\varphi(z, s_1, \tilde s')=0,$ 
$\alpha, \beta \in \Omega^{n-1}$ in ${\mathcal B}_\varphi(\tilde s')$, entails the existence of $\alpha' \in \Omega^{n-1}$
such that $d\omega =d\alpha' \wedge d\varphi(z, s_1, \tilde s'), $ this in turn together with the de Rham lemma entails
$ \omega =\alpha' \wedge d\varphi(z, s_1, \tilde s') + d\beta'$
for some  $\beta' \in \Omega^{n-1}$

To see the surjectivity, it is enough to check that for every $\gamma \in \Omega^{n}$ the equation $d\omega =\gamma$
is solvable.
{\bf Q.E.D.}

Let us introduce a module for $\tilde s'=(\tilde s_2, \cdots, \tilde s_m)\in U$,
$$Q_\varphi( \tilde s')  := \frac{\Omega_{\bC^n}^{n}}{d_z\varphi(z, s_1, \tilde s')
\wedge  \Omega_{\bC^n}^{n-1}} \cong \frac{{\bC}[z]}{(d_z\varphi(z, s_1, \tilde s')){\bC}[z]}, \eqno(2.7)$$
 that is a  free $\bC[s_1]$ module of rank $\mu$  because it is isomorphic to
$$ \oplus_{\{s_1: Z_{( s_1, \tilde s')}\; \rm{singular} \}} \oplus_{z: singular\; points\; on \;Z_{( s_1, \tilde s')}} \bC^{\mu(z)},$$
with $\mu(z) :$ the Milnor number of the singular point $z \in Z_{( s_1, \tilde s')}$.  Let us denote its basis by 
 $$ \{g_1 dz, \cdots, g_\mu dz \},  \eqno(2.8)$$
such that the polynomials  $\{g_1 (z), \cdots, g_\mu (z) \}$ 
consist a basis of the RHS of $(2.7)$ as a free $\bC[s_1]$ module.

According to \cite[p.218, lines 5-6]{Br88}
the following is a locally trivial fibration,
$$ Z_{(s_1, s')} \rightarrow (s_1, s') \in \bC \times U \setminus  D_\varphi,$$ by the definition $(2.5),$ $(2.6).$This yields the next statement.
\begin{corollary}
We can choose a basis $\{\omega_1, \cdots, \omega_\mu \}$ of ${\mathcal P}_\varphi(\tilde s')$ independent of
 $\tilde s' \in U.$
\end{corollary}

Due to the construction of $U$, we can consider the ring ${\mathcal O}_U$ of holomorphic functions on $U.$
By the analytic continuation with respect to the parameter $s' \in U,$ we see the following.

\begin{lemma} The modules
${\mathcal B}_\varphi(s')$,  ${\mathcal P}_\varphi(s')$, $Q_\varphi(s')$ are free 
 $\bC[s_1] \otimes {\mathcal O}_U$ modules of rank $\mu$. 
 \label{lemma3}
\end{lemma}

As the deformation polynomials $e_{1}, \cdots, e_{\mu}$ arise from the special form of $\Psi(x,t,z)$ we are obliged to impose the following assumption.

\vspace{1pc}
\noindent
{\bf Assumption II}
 We assume that  we can adopt $e_i(z)$ of $(1.5),$ $(1.6)$ as $g_i(z)$ in $(2.8)$ 
$i=1,\cdots,\mu$ and they serve as a basis of  $Q_\varphi(s')$ as a free 
 $\bC[s_1] \otimes {\mathcal O}_U$ module. 

\vspace{1pc}
For the sake of simplicity, let us denote by  $mod (d_z(\varphi(z,0) +\sum_{j=2}^m s_j e_j(z)))$ the residue class modulo
the ideal $(d_z(\varphi(z,0) +\sum_{j=2}^m s_j e_j(z))) \bC[z,s_1] \otimes {\mathcal O}_U$ in $\bC[z,s_1] \otimes {\mathcal O}_U$.  By virtue of the freeness of $Q_\varphi(s')$,
this residue class is uniquely determined.
Our assumption  $(1.5), (1.6)$ together with the Weierstrass preparation theorem gives us a decomposition as follows,
$$ (\varphi(z,0) +\sum_{j=2}^m s_j e_j(z))\cdot \frac{\partial \varphi(z,s) }{\partial s_i} $$
$$\equiv  \sum_{\ell=1}^\mu \sigma_{i}^\ell(s') \frac{\partial \varphi(z,s) }{\partial s_\ell} \;\; mod (d_z(\varphi(z,0) +\sum_{j=2}^m s_j e_j(z))), \; 1 \leq i \leq \mu \eqno(2.9)$$
$$  \frac{\partial \varphi(z,s) }{\partial s_i} \equiv 
\sum_{\ell=1}^\mu \sigma_{i}^\ell(s') \frac{\partial \varphi(z,s) }{\partial s_\ell} \;\; mod (d_z(\varphi(z,0) +\sum_{j=2}^m s_j e_j(z))),$$ 
$$\; \mu+1 \leq i \leq m,  \eqno(2.10)$$ 
with $\sigma_{i}^\ell(s') \in {\mathcal O}_U.$
In fact, according to an argument used in ~\cite[Theorem A4]{bru}, ~\cite[Proposition 2]{Tan07} (both treat liftable vector fields in local case but they apply to our situation), the following vector fields are tangent to the discriminant $D_\varphi,$
$$ \vec v_i:=(s_1+\sigma_{i}^i (s') )\frac{\partial  }{\partial s_i}+ \sum_{\ell=1, \ell \not =i}^\mu \sigma_{i}^\ell(s') \frac{\partial \varphi(z,s) }{\partial s_\ell}, \; 1 \leq i \leq \mu \eqno(2.11)$$
Here we recall the Assumption I, (i) that allows us to adopt $(s_1, s_2, \cdots, s_\nu),$ $\nu \geq \mu$
as the local coordinates of $\bC \times U$.
$$   \vec v_i:=-\frac{\partial  }{\partial s_i}+ \sum_{\ell=1}^\mu \sigma_{i}^\ell(s') \frac{\partial }{\partial s_\ell}, \;\; \mu+1 \leq i \leq \nu,  \eqno(2.12) $$ 

We remark here that the importance of the liftable vector field in the studies of $A_\mu$ singularity discriminant has been pointed out by  ~\cite[\S 3]{Ar76}.

Evidently they are linearly independent over $\bC[s_1]\otimes {\mathcal O}_U$ because of the presence of the term $s_1\frac{\partial  }{\partial s_i}$
for every $1 \leq i \leq \mu $ and $-\frac{\partial  }{\partial s_i}$ for $\mu+1 \leq i \leq \nu$. Therefore they  form a $\bC[s_1]\otimes {\mathcal O}_U$free module of rank $\nu$.
Let us introduce the following matrix of which the $i-$th row corresponds to the vector $\vec v_i$.
$$ \Sigma(s):=
\left (\begin {array}{ccccccccccccc} s_1+\sigma_{1}^1 (s') & \sigma_{1}^2 (s')& \cdots &\sigma_{1}^\mu(s') &0&\cdots&0&0
\\\noalign{\medskip}\sigma_{2}^1 (s')&s_1+\sigma_{2}^2 (s')&\cdots & \sigma_{2}^\mu(s')&0&\cdots&0&0
\\\noalign{\medskip}\vdots &\vdots&\ddots&\vdots&\vdots&\vdots&\vdots&\vdots
\\\noalign{\medskip}\sigma_{\mu}^1 (s')&\sigma_{\mu}^2 (s')&\cdots & s_1+\sigma_{\mu}^\mu(s')&0&\cdots&0&0
\\\noalign{\medskip}\sigma_{\mu+1}^1 (s')&\sigma_{\mu+1}^2 (s')&\cdots & \sigma_{\mu+1}^\mu(s')&-1&\cdots&0&0
\\\noalign{\medskip}\vdots&\vdots&\cdots & \vdots&\vdots&\ddots&\vdots&\vdots
\\\noalign{\medskip}\sigma_{\nu-1}^1 (s')&\sigma_{\nu-1}^2 (s')&\cdots & \sigma_{\nu-1}^\mu(s')&0&\cdots&-1&0
\\\noalign{\medskip}\sigma_{\nu}^1 (s')&\sigma_{\nu}^2 (s')&\cdots & \sigma_{\nu}^\mu(s')&0&\cdots&0&-1
\end {array}\right ).   \eqno(2.13)$$
In fact the following $\mu \times \mu$ submatrix of $\Sigma(s)$ contains the essential geometrical informations on
$D_\varphi.$ 
 $$ \tilde \Sigma(s):=
\left (\begin {array}{ccccccccccccc} s_1+\sigma_{1}^1 (s') & \sigma_{1}^2 (s')& \cdots &\sigma_{1}^\mu(s')
\\\noalign{\medskip}\sigma_{2}^1 (s')&s_1+\sigma_{2}^2 (s')&\cdots & \sigma_{2}^\mu(s')
\\\noalign{\medskip}\vdots &\vdots&\ddots&\vdots
\\\noalign{\medskip}\sigma_{\mu}^1 (s')&\sigma_{\mu}^2 (s')&\cdots & s_1+\sigma_{\mu}^\mu(s')
\end {array}\right ).   \eqno(2.14)$$

\begin{theorem}
1) The algebra $Der_{\bC \times U}(log\; D_\varphi)$ 
of tangent fields to $D_\varphi$ as a free  $\bC[s_1]\otimes {\mathcal O}_U$ is generated by the vectors $v_i,$ $1 \leq i \leq \nu$ of $(2.11),$ $(2.12).$ 

2) 
The discriminantal loci $ D_\varphi$ is given by 
the equation 
$$D_\varphi=\{s \in  \bC \times U : det\; \tilde \Sigma(s)=0 \}.$$

3) The preimage of $D_\varphi$ by the mapping $\iota$
contains the wave front $BW=\cup_{t \in \bC} W_t \subset \bC^{n+1}$ i.e. $BW \subset \iota^{-1}(D_\varphi).$ 
\label{theorem2}
\end{theorem}

{\bf Proof} The tangency of vector fields $\vec v_i$'s to   $D_\varphi$  and their independence over $\bC[s_1]\otimes {\mathcal O}_U$ have already been shown.
 
  We shall follow the argument by \cite[Theorem 3.1]{Gory1}. First we shall prove $2).$ By virtue of the tangency of  $\vec v_i$'s to   $D_\varphi$ and the equality,
 $$  \vec v_1 \wedge \cdots\wedge \vec v_\nu = det \;\Sigma(s) {\partial_{s_1}} \wedge \cdots \wedge {\partial_{s_\nu}},$$
 the function $det \;\Sigma(s)$ shall vanish on $D_\varphi$. The statement on $ Q_\varphi(s')$ of the Lemma ~\ref{lemma3} 
 tells us that
 $$ \sharp \{s \in \bC \times U: s_1=const\;\cap D_\varphi\} =\mu, $$
 in taking the multiplicity into account.
 
 From $(2.13),$ $(2.14)$ we see that 
 $$  \pm det \;\Sigma(s) = det \;\tilde \Sigma(s)= s_1^\mu + d_1(s') s_1^{\mu-1} + \cdots + d_\mu(s'),$$
with $d_i(s') \in  {\mathcal O}_U,$ $1 \leq i \leq \mu$. Thus the determinant $det \;\tilde \Sigma(s)$ that turns out to be a Weierstrass polynomial in $s_1$,
shall be divided by the defining equation of  $D_\varphi$ which turns out to be also a Weierstrass polynomial in $s_1$ of degree $\mu$.
This proves $2).$

Now we shall show that every vector $\vec v$ tangent to $D_\varphi$ admits a decomposition like
$$ \vec v = \sum^\nu_{i=1} a_i(s) \vec v_i,$$
for some $a_i(s) \in \bC[s_1]\otimes {\mathcal O}_U.$ For every $i$
the following expression shall vanish on $D_\varphi$, because of the tangency of
all vectors taking part in it,
$$  \vec v_1 \wedge \cdots \wedge \vec v_{i-1} \wedge \vec v \wedge \vec v_{i+1} \wedge \cdots \wedge \vec v_\nu.$$
Therefore there exists  $a_i(s) \in \bC[s_1]\otimes {\mathcal O}_U$ such that the above expression equals to
$a_i(s)  det \Sigma(s) {\partial_{s_1}} \wedge \cdots \wedge {\partial_{s_m}}.$
This means that the vector $\vec v$ $-$ $\sum^\nu_{i=1}$ $a_i(s)$ $\vec v_i$  defines a zero vector at every $s \not \in D_\varphi,$
as the vectors $\vec v_1, \cdots, \vec v_\nu $ form a frame outside  $D_\varphi$. By the continuity argument on holomorphic functions, we see that the decomposition
holds everywhere on $ \bC \times U.$

The statement $3)$ follows from Lemma ~\ref{lemma1}, $(1.4)$ and the definition $(1.7)$ of  $\iota.$
{\bf Q.E.D.}

{
\center{\section {Gauss-Manin system for a tame polynomial}}
}
In this section,
we willl show that 
the above matrix $\tilde \Sigma(s)$, $(2.14)$ can be obtained as the coefficient of the Gauss-Manin system defined for a tame polynomial
$\varphi(z,s)$.

According to Lemma ~\ref{lemma3}, every $\omega \in  {\mathcal P}_\varphi(s')$ admits a unique decomposition as follows,
$$ \omega = \sum^\mu_{i=1} a_i(s) \omega_i, \;\;\;\; s \in \bC \times U. \eqno(3.1)$$
A generalisation of Theorem 0.2 of \cite{G} tells us that the following equivalence holds for every holomorphic $n-1$
form $\omega,$
$$ \forall s \in \bC\times U, \omega|_{Z_{s}} =0 \;\; {\rm in}\;\; H^{n-1}(Z_{s}) \Leftrightarrow  \omega=0 \;\; \rm{in}\;\; {\mathcal P}_\varphi(s').\eqno(3.2)$$
The above statement $(3.2)$ for every $n \geq 2$ was given by Corollary 10.2 of \cite{Sab06} after an argument quite different from that in  \cite{G} \S 2.

This theorem yields a corollary that
ensures us the following equality for every vanishing cycle $\delta(s) \in H_{n-1}(Z_{s}),$
$$ \int_{\delta(s)} \omega = \sum^\mu_{i=1} a_i(s)  \int_{\delta(s)}\omega_i, s \in \bC \times U, \eqno(3.3)$$
for some $a_i(s) \in \bC[s_1]\otimes {\mathcal O}_U$, $1 \leq i \leq \mu$.
To show this along with the argument by \cite{G}, we simply need to replace his Lemma 2.2 by \cite[Corollary 0.7]{DiS}.

Here we remark that for the basis of $\{e_1(z)dz, \cdots, e_\mu(z)dz \}$ of $Q_\varphi(\tilde s')$
we can choose the basis  $\{\omega_1, \cdots, \omega_\mu \}$ of ${\mathcal P}_\varphi(\tilde s')$
such that 
$$ d\omega_i = e_i(z) dz +d_z \varphi(z,s) \wedge \epsilon_i, \eqno(3.4)$$ 
for some $\epsilon_i  \in \Omega^{n-1}.$ That is to say, for every $\omega \in  \Omega^{n-1} $
we can find the following two types of decomposition
$$ \omega = \sum _{i=1}^\mu c_i(s') d\omega_i  + d_z \varphi(z,s) \wedge d\xi,$$
$$  = \sum _{i=1}^\mu c_i(s') (e_i(z) dz +d_z \varphi(z,s) \wedge \epsilon_i)  + d_z \varphi(z,s) \wedge \eta, \eqno(3.5) $$
for some
$c_i(s') \in  {\mathcal O}_U,$
$\xi \in \Omega^{n-2} \otimes {\mathcal O}_U$, $\eta \in \Omega^{n-1} \otimes {\mathcal O}_U.$
This is a reformulation of Lemma ~\ref{lemma3}.

As E.Brieskorn ~\cite{Br70} showed,  the following equality holds if we understand it as the property of
the holomorphic sections in the cohomology bundle $H^{n-1}(Z_s)$ defined as the Leray's residue $\omega/d_z \varphi(z,s)$
for $\omega \in \Omega^n,$
$$ (\frac {\partial }{\partial s_1})^{-1} d\eta =  d_z \varphi(z,s) \wedge \eta.  $$
This yields  that
$$(\frac {\partial }{\partial s_1})^{-1} {\mathcal B}_\varphi(\tilde s')= d_z \varphi(z,s) \wedge \Omega^{n-1}/ d_z \varphi(z,s) \wedge d\Omega^{n-2},$$
$$  Q_\varphi(\tilde s')= {\mathcal B}_\varphi(\tilde s')/ (\frac {\partial }{\partial s_1})^{-1} {\mathcal B}_\varphi(\tilde s'),$$
and we see that $\{e_1(z)dz, \cdots, e_\mu(z)dz \}$ is a basis of ${\mathcal B}_\varphi(\tilde s')$ as an $ {\mathcal O}_U[(\frac {\partial }{\partial s_1})^{-1}] $
module.

For $\omega_i$'s chosen in $(3.4)$ we have a  decomposition in $Q_\varphi(\tilde s')$ as follows,
$$ (\varphi(z,s) -s_1) d\omega_i = \sum^\mu_{\ell=1} \sigma^\ell_i(s') d\omega_\ell + d_z \varphi(z,s) \wedge \eta_i, \;\;\;\;1 \leq i \leq \mu \eqno(3.6)$$
$ \eta_i \in \Omega^{n-1} \otimes {\mathcal O}_U.$  We see that $(3.6)$ is equivalent to $(2.9)$ in view of $(3.5).$
This relation immediately entails the following equality for every  $\delta(s) \in H_{n-1}(Z_{s}),$
$$  s_1 \frac{\partial }{\partial s_1}\int_{\delta(s)}\omega_i + \sum^\mu_{\ell=1} \sigma^\ell_i(s') \frac{\partial }{\partial s_1}\int_{\delta(s)}\omega_\ell + \int_{\delta(s)} \eta_i =0,\eqno(3.7)$$
in view of the fact $ \int_{\delta(s)} \varphi(z,s)  \frac{d\omega_i}{d_z \varphi(z,s)} =0$ and the Leray's residue theorem
$$  \frac{\partial }{\partial s_1}\int_{\delta(s)}\omega_i=  \int_{\delta(s)} \frac{d\omega_i}{d_z \varphi(z,s)}.$$

After $(3.3)$, every $\int_{\delta(s)} \eta_i$ admits an unique  decomposition 
$$ \int_{\delta(s)} \eta_i = \sum^\mu_{j=1} b_i^j(s)  \int_{\delta(s)}\omega_j,  s \in \bC \times U, \eqno(3.8)$$
for some $b_i^j(s) \in \bC[s_1]\times {\mathcal O}_U$,  $1 \leq i,j \leq \mu$.

Let us  consider a vector of fibre
integrals 
$$\bI_{Q}:= ^t(\int_{\delta(s)}\omega_1, \cdots, \int_{\delta(s)}\omega_\mu). \eqno(3.9) $$
In summary we get
\begin{theorem}
1) For a vector $\bI_{Q}$, $(3.7)$ we have the following Gauss-Manin system 
$$ \tilde \Sigma \cdot \frac{\partial }{\partial s_1} \bI_{Q}+ B(s) \bI_{Q}=0, \eqno(3.10)$$
where $B(s) = \left( b_i^j(s) \right)_{1 \leq i,j \leq \mu}$ for functions determined in $(3.8).$

2) The discriminantal loci $D_\varphi$ of the  tame polynomial $\varphi(z,s) ,$ $s \in \bC \times U$
has an expression,
$$D_\varphi=\{s \in  \bC \times U : det\; \tilde \Sigma(s)=0 \},$$
that corresponds to the singular loci of the system $(3.10)$.
\label{theorem3} \end{theorem}

\begin{remark}
To see that the two statements on $D_\varphi$ do not mean a simple coincidence, one may consult ~\cite[Theorem 2.3]{Tan07}
and find a description of the Gauss-Manin system for Leray's residues by means of the tangent vector fields to
the discriminant loci.
\end{remark}

{\section {Free and almost free wave fronts}}

Now we recall that the freeness of
$Der_{\bC \times U}(log\; D_\varphi)$ as a $\bC[s_1]\otimes {\mathcal O}_U$ module, proven in the Theorem ~\ref{theorem2},
means that  $D_\varphi$ defines a free divisor  (in the sense of K.Saito) in the neighbourhood of 
every point $s \in   D_\varphi.$ We define the logarithmic tangent space $T^{log}_{s} D_\varphi$ to $ D_\varphi$ at $s$:
$$ T^{log}_{s} D_\varphi = \{\vec v(s) : \vec v(s) \in Der_{\bC \times U}(log\; D_\varphi)_s\} \eqno(4.1)$$
We follow the presentation  on the free and almost free divisors by D.Mond ~\cite{Mond} and J.N.Damon ~\cite{Damon}.
To discuss when the big wave front $BW$ becomes a free divisor, we need to make use of the notion of algebraic transversaliy.
We recall here the Assumption I, (ii) on the image of the mapping $\iota$ that entails the following inclusion relation,
$$ d_{x,t}\iota (T_{(x,t)} \bC^{n+2}) \subset  T_{\iota(x,t)}(\bC \times U),$$
for $(x,t) $ in the neighbourhood of $(x_0,t_0).$
\begin{definition} 
The mapping $\iota$ is algebraically transverse to $D_\varphi$ at $(x_0,t_0) \in \bC^{n+2}$
if and only if  
$$ d_{x,t}\iota (T_{(x_0,t_0)} \bC^{n+2}) + T^{log}_{\iota(x_0,t_0)} D_\varphi =  T_{\iota(x_0,t_0)}(\bC \times U). \eqno(4.2)$$
\end{definition}
\begin{lemma}(~\cite{Mond}  Jacobian criterion for freeness)
The divisor $\iota^{-1}(D_\varphi)$ is free if and only if $\iota$ is algebraically transverse to $ D_\varphi$.
\label{lemma41}
\end{lemma}
To state a criterion of the freeness of $\iota^{-1}(D_\varphi)$, we need the following  $m \times (\nu+n+2)$  matrix $T(x,t).$
 $$T(x,t)=$$ {\tiny $${
\left (\begin {array}{ccccccccccccc} s_1+\sigma_{1}^1 (s'(x,t)) & \cdots &\sigma_{1}^\mu(s'(x,t)) &0&\cdots&0&\cdots&0
\\\noalign{\medskip}\sigma_{2}^1 (s'(x,t))&\cdots & \sigma_{2}^\mu(s'(x,t))&0&\cdots&0&\cdots&0
\\\noalign{\medskip}\vdots &\ddots&\vdots&\vdots&\vdots&\vdots&\cdots&\vdots
\\\noalign{\medskip}\sigma_{\mu}^1 (s'(x,t))&\cdots & s_1+\sigma_{\mu}^\mu(s'(x,t))&0&\cdots&0&\cdots&0
\\\noalign{\medskip}\sigma_{\mu+1}^1 (s'(x,t))&\cdots & \sigma_{\mu+1}^\mu(s'(x,t))&-1&\cdots&0&\cdots&0
\\\noalign{\medskip}\vdots&\cdots & \vdots&\vdots&\ddots&\vdots&\cdots&\vdots
\\\noalign{\medskip}\sigma_{\nu}^1 (s'(x,t))&\cdots & \sigma_{\nu}^\mu(s'(x,t))&0&\cdots&-1&\cdots&0
\\\noalign{\medskip}s_1(x,t)_{x_1}&\cdots &s_\mu(x,t)_{x_1}&s_{\mu+1}(x,t)_{x_1}&\cdots&s_{\nu}(x,t)_{x_1}&\cdots&s_{m}(x,t)_{x_1}
\\\noalign{\medskip}s_1(x,t)_{x_2}&\cdots &s_\mu(x,t)_{x_2}&s_{\mu+1}(x,t)_{x_2}&\cdots&s_{\nu}(x,t)_{x_2}&\cdots&s_{m}(x,t)_{x_2}
\\\noalign{\medskip}\vdots &\vdots&\vdots&\vdots&\vdots&\vdots&\vdots&\vdots
\\\noalign{\medskip}s_1(x,t)_{x_{n+1}}&\cdots &s_\mu(x,t)_{x_{n+1}}&s_{\mu+1}(x,t)_{x_{n+1}}&\cdots&s_{\nu}(x,t)_{x_{n+1}}&\cdots&s_{m}(x,t)_{x_{n+1}}
\\\noalign{\medskip}s_1(x,t)_{t}&\cdots &s_\mu(x,t)_{t}&s_{\mu+1}(x,t)_{t}&\cdots&s_{\nu}(x,t)_{t}&\cdots&s_{m}(x,t)_{t}
\end {array}\right )}.$$} $$\eqno(4.3)$$
The first $\nu$ rows of the $T(x,t)$ correspond to those of $\Sigma(\iota(x,t))$ while the $(\nu+i)-$th row corresponds to
$\frac{\partial}{\partial x_i} \iota(x,t),$ $1 \leq i \leq n+1$ and the last row to $\frac{\partial}{\partial t}\iota(x,t)$ for $\iota(x,t)$ of $(1.7).$

The Lemma ~\ref{lemma41} yields immediately the following statement in view of the Theorem ~\ref{theorem2}.
\begin{proposition}
The divisor $\iota^{-1}(D_\varphi)$ is free in the neighbourhood of  $(x,t)$ if and only if
 $rank \; T(x,t) \geq \nu.$ 
\label{proposition41}
\end{proposition}

\begin{remark}
It is well known that the discriminant of a $\mathcal K-$ versal deformation of a hypersurface singularity defines a free divisor (K. Saito, E. Looijenga).
Therefore, the cases for which the $\mathcal K-$ versality has been proven in ~\cite{Hase}, $n=2$ give rise to  free wave fronts.  
\end{remark}

After Theorem ~\ref{theorem2}, in the neighbourhood of each of its point $s$, the  hypersurface $ D_\varphi$ 
defines a germ of free divisor. 
\begin{definition} 
The germ of hypersurface  $\iota^{-1}(D_\varphi)$ at $(x_0,t_0) \in \bC^{n+2}$
is an almost free divisor based on the germ of free divisor $D_\varphi$ at $\iota(x_0,t_0) \in \bC \times U$
if there is a map $i_0: \iota^{-1}(D_\varphi) \rightarrow D_\varphi$ which is algebraically transverse to $D_\varphi$
except at $(x_0,t_0)$ such that  $\iota^{-1}(D_\varphi) = i_0^{-1}( D_\varphi).$
\end{definition}

In view of this definition, we get a  criterion so that $\iota^{-1}(D_\varphi)$ be an almost free divisor.
\begin{proposition}
The germ of hypersurface  $\iota^{-1}(D_\varphi)$ at $(x_0,t_0) \in \bC^{n+2}$
is an almost free divisor based on the germ of free divisor $D_\varphi$ at $\iota(x_0,t_0) \in \bC \times U$ if 
the following inequality holds at an isolated point $(x_0,t_0) \in \iota^{-1}(D_\varphi)$,
$$rank\;T(x_0,t_0) < \nu,  \eqno(4.4)$$
while at other points $(x,t) \not = (x_0,t_0)$ in the neighbourhood of $(x_0,t_0)$ , the inequality $rank\; T$ $ (x, t)$ $\geq$ $\nu$ holds.
\label{proposition42}
\end{proposition}

As we shall see in the Example 5.1 below, it is quite difficult to verify that the condition $(4.4)$ is satisfied at an isolated point. We can give a sufficient condition on the violation of algebraic transversality condition at an isolated point as follows.

\begin{proposition} Assume that $(4.4)$ holds at the focal point $(x_0,t_0).$
For a mapping $\iota$ with $rank\;d_{x,t} \iota(x_0,t_0)= n+1$, if the following inequality $(4.5)$ is satisfied only for  $(\xi,\tau)=$ $(0,0)$, $(x_0,t_0)$is an isolated algebraically non-transverasal point. 
$$  T(x_0,t_0) +\tau  \frac{\partial T}{\partial t}(x_0,t_0)+\sum_{j=1}^n\xi_j \frac{\partial T}{\partial x_j}(x_0,t_0) < \nu. \eqno(4.5)$$
\label{proposition44}
\end{proposition}

The proof follows directly from the lower semi-continuous property of the  $rank\;T$ $(x,t).$

{\section {Examples}}

{\bf 1. Wave propagation on the plane}

 Let us consider the following initial wave front on the plane
$Y:=\{(z,u) \in \bC^{2};a z^2 +z^4+u=0  \}, z=$ i.e. $F (z)= a z^2 +z^4 $ for some real non-zero constant
$a.$ In this case our phase function has the following expression
$$\Psi(x,t,z)=  (x_1 + a z^2 + z^4 + (x_2 - z) (2 a z + 4 z^3))^2 - 
 t^2 (1 + (2 a z + 4 z^3)^2),$$
$$ =-t^2 + x_2^2 + 4 a x_1 x_2 z + (-4 a^2 t^2 + 4 a^2 x_1^2 - 2 a x_2) z^2$$ $$(-4 a^2 x_1 + 8 x_1 x_2) z^3 + (a^2 - 16 a t^2 + 16 a x_1^2 - 6 x_2) z^4 $$ $$- 
 20 a x_1 z^5 + (6 a - 16 t^2 + 16 x_1^2) z^6 - 24 x_1 z^7 + 9 z^8. \eqno(5.1)$$
It is easy to see that $(x_1,x_2,t) =(0,-\frac{1}{2a}, \frac{1}{2a})$ 
is a focal point with  a singular point $(z,u)=(0,0)$
and the Milnor number $\mu(0)= 3$  ($A_3$ singularity i.e. the swallow tail) if 
$ a \not = 1$ and $\mu(0)= 5$ ($A_5$ singularity) if $a=1$,
$$\Psi(0,-\frac{1}{2a}, \frac{1}{2a},z) =(-(1/a) + a^2) z^4 + (-(4/a^2) + 6 a) z^6 + 9 z^8.  \eqno(5.2)$$
 
 The quotient ring $(1.5)$ for this $\Psi(0,-\frac{1}{2a}, \frac{1}{2a},z)$ has dimension $\mu=7$.

Especially we can choose   $e_i = z^{i-1},$ $i=1,\cdots, 7$  as the basis of $(2.8).$
Now, in view of $(5.1)$ we introduce an additional deformation polynomial
$e_8=z^7,$
together with entries of the mapping $\iota$ $(1.7),$
$$ s_1=-t^2 + x_2^2, s_2= 4 a x_1 x_2, s_3=-4 a^2 t^2 + 4 a^2 x_1^2 - 2 a x_2, s_4=-4 a^2 x_1 + 8 x_1 x_2,$$ $$ s_5=a^2 - 16 a t^2 + 16 a x_1^2 - 6 x_2, s_6=-20 a x_1, 
s_7=6 a - 16 t^2 + 16 x_1^2, s_8= -24 x_1. \eqno(5.3)$$
$$\varphi(z,s) = 9 z^8 + \sum^8_{i=1} s_iz^{i-1}.$$
In this case, the constructible set $U$ of the Assumption I,(i) coincides with $\bC^7.$

 At the focal point  $(x,t)= (0,-\frac{1}{2a}, \frac{1}{2a})$ the matrix   $\iota^\ast(\Sigma)(0,-\frac{1}{2a}, \frac{1}{2a})$  has the following form with rank 5 if $a \not=1$ and rank 3 if $a=1$.
{\tiny $$ \left (\begin {array}{ccccccccccccc}0& 0& 0& 0& (-1 + a^3)/(2 a)& 0& -(1/a^2) + (3 a)/2& 0
\\\noalign{\medskip}0& 0& 0& A_1& 0& A_2 & 0& 0
\\\noalign{\medskip}0& 0& 0& 0& A_1& 0& A_2& 0
\\\noalign{\medskip}0& 0& 0& A_3& 0&A_4& 0& 0
\\\noalign{\medskip}0& 0& 0& 0& A_3& 0& A_4& 0
\\\noalign{\medskip}0& 0& 0& A_5& 0&A_6& 0& 0 
\\\noalign{\medskip}0& 0& 0& 0& A_5& 0& A_6& 0 
\\\noalign{\medskip}0& 0& 0& 4 (-1 + a^3)/a& 0& 6 (-(4/a^2) + 6 a)& 0& 72\end {array} \right)
$$}
  $$ \eqno(5.4)$$
  where  
  $A_1= \frac{-(2 - 5 a^3 + 3 a^6)}{36 a^3},$
 $A_2= \frac{-(4 - 6 a^3 + 3 a^6)}{12 a^4},$
 $A_3= \frac{-4 + 10 a^3 - 9 a^6 + 3 a^9}{216 a^5},$
$A_4= \frac{(-2 + a^3)^2 (-2 + 3 a^3)}{72 a^6},$
$A_5=-\frac{(-2 + a^3)^2 (2 - 5 a^3 + 3 a^6)}{1296 a^7},$
$A_6=-\frac{16 - 56 a^3 + 68 a^6 - 30 a^9 + 3 a^{12}}{432 a^8}.$
We see that $A_1=A_3=A_5=0$ for $a=1.$

Thus together with the data    
     $$d_{x,t} \iota(0,-1/2a, 1/2a) \eqno(5.5)$$
$$=\left (\begin {array}{ccccccccccccc}0& -2& 0& -(4/a) - 4 a^2& 0& -20 a& 0& -24
\\\noalign{\medskip}0&	-(1/a)& 0& -2 a& 0& -6& 0& 0
\\\noalign{\medskip} -(1/a)& 0& -4 a& 0& -16& 0& -(16/a)& 0
\end {array} \right)
$$
we conclude that $rank\; T(0,-\frac{1}{2a}, \frac{1}{2a})= 8 =\nu$ if $a \not =1$. Therefore after Proposition ~\ref{proposition41}, the germ of the big wave front $BW$ defines a free divisor in the neighbourhood of the focal point $(0,-1/2a, 1/2a)$ for $a  \not=1.$

In the case $a =1,$ $rank \; \iota^\ast(\Sigma )(0,-1/2, 1/2)=$ $rank \; \iota^\ast(\tilde \Sigma )(0,-1/2, 1/2)+1 =3$
and $$rank\; T(0,-1/2, 1/2)= 6 <8.  \eqno(5.6)$$
That is to say the mapping $\iota$ is not algebraically transverse at the focal point $(0,-1/2, 1/2).$
 Now we shall see that the focal point $(0,-1/2, 1/2)$
is an isolated point after the following reasoning.

At first we remark that $(5.4)$, $(5.5)$ entail the following relation.
$$ span_{\bC}\{ v_1(\iota(0,-1/2, 1/2)), \cdots, v_8(\iota(0,-1/2, 1/2)) \}$$ $$\cap span_{\bC}\{  \frac{\partial \iota}{\partial t}, \frac{\partial \iota}{\partial x_1}, \frac{\partial \iota}{\partial x_2}\}_{(0,-1/2, 1/2)} =\{0\}.  $$
This means that the integral variety germ of the vectof fields $\{ v_1(s),$ $ \cdots,$ $ v_8(s) \}$(i.e. the $A_5$ singularity stratum of $D_{\varphi,\iota(0,-1/2, 1/2)}$) and the image $\iota (\bC^3,$ $(0,$ $-1/2,$ $1/2))$ intersect transversally (in the usual sense) at the point $\iota(0,-1/2, 1/2).$ This does not ensure the isolation of the algebraic non-transversality.
We still need to show that no $A_4$ singularity stratum on the wave front is adjacent to the focal point.

\vspace{2pc} 

{\bf First proof of the isolated property of the focal point }

If $\Psi(x,t,z)$ had a $A_4$ singularity point 
in the neighbourhood of the focal point $(0,-1/2, 1/2)$ with 
$A_5$ singularity, the rank of  $T(x,t)$ would be 7 ($<8$) there. At such a $A_4$ singularity point, the condition of
Proposition~\ref{proposition42} is not satisfied. In this situation, the algebraic transversality would be violated on a non-discrete set adjacent to the focal point. If we show that the $A_4$ singularity stratum of $D_\varphi$
is not contained in the image $\iota($ $\bC^3,$ $(0,$ $-1/2,$ $1/2)$ $) \setminus \iota(0,-1/2, 1/2)$ $\subset $ $\bC \times U$, it would mean that  no  $A_4$ singularity appears on the wave front. Consequently it proves that $rank T(x,t) \geq 8$ for every $(x,t) \not = (0,-1/2, 1/2).$

A deformation of the polynomial $(5.2)$ with
$A_k$ ($k \geq 4$) singularity near the origin can be given by 
$$(z + w_1)^5(q_1 + q_2z + q_3z^2 + 9z^3),$$
for $(w_1,q_1,q_2,q_3) \approx (0,0,2,0).$
In other words, the union of $A_k$ singularity ($k \geq 4$) strata in $\bC \times U$
has the following 4-parameter representation,
$$ (s_1,\cdots, s_8)=$$
{\small $ (q_1 w_1^5, 5 q_1 w_1^4 + q_2 w_1^5, 10 q_1 w_1^3 + 5 q_2 w_1^4 + q_3 w_1^5, 
 10 q_1 w_1^2 + 10 q_2 w_1^3 + 5 q_3 w_1^4 + 9 w_1^5,$ }
{\small  $5 q_1 w_1 + 10 q_2 w_1^2 + 10 q_3 w_1^3 + 45 w_1^4, 
 q_1 + 5 q_2 w_1 + 10 q_3 w_1^2 + 90 w_1^3, q_2 + 5 q_3 w_1 + 90 w_1^2, 
 q_3 + 45 w_1).$}

Four vectors below span the tangent to this union set at the point $(w_1,q_1,$ $q_2,$ $q_3)$ $=$ $(0,$ $0,2,0)$,
$$(0, 0, 0, 0, 0, 1, 0, 0),(0, 0, 0, 0, 0, 0, 1, 0),(0, 0, 0, 0, 0, 0, 0, 1),$$ $$(0, 0, 0, 0, 0, 10, 0, 45).$$
This three dimensional tangent space and $span_{\bC}\{  \frac{\partial \iota}{\partial t}, \frac{\partial \iota}{\partial x_1}, \frac{\partial \iota}{\partial x_2}\}_{(0,-1/2, 1/2)}$ have a common vector subspace $\{0\}.$
This yields that the union of $A_k$ singularity ($k \geq 4$) strata and $\iota(\bC^3,(0,-1/2, 1/2))$ intersect only at the point
$\iota(0,-1/2, 1/2).$
Hence  $A_k$ ($k \geq 4$) does not appear in $\iota(\bC^3,(0,-1/2, 1/2))$ $ \setminus $ $ \iota ( $ $0,$ $-1/2,$ $ 1/2$ $).$


\vspace{2pc}

{\bf Second proof of the isolated property of the focal point}

We verify the isolated algebraically non-transversal property at $(0,$ $-1/2,$ $1/2)$ by means of Proposition~\ref{proposition44}.
In this case the LHS of $(4.5)$ becomes
$$T(0,-\frac{1}{2}, \frac{1}{2})+\tau  \frac{\partial T}{\partial t}(0,-\frac{1}{2}, \frac{1}{2})+\xi_1 \frac{\partial T}{\partial x_1}(0,-\frac{1}{2}, \frac{1}{2})+\xi_2 \frac{\partial T}{\partial x_2}(0,-\frac{1}{2}, \frac{1}{2}),\eqno(5.6)$$
Among the above matrices $T(0,-\frac{1}{2}, \frac{1}{2})$is already given by $(5.4),(5.5).$
Other derivatives are calculated as follows,
 
\newpage
$$  \frac{\partial T}{\partial t}(0,-1/2, 1/2)=$$
{\tiny $$  \left (\begin{array}{ccccccccccccc}-1& 0& -3& 0& -8& 0& -4& 0
 \\\noalign{\medskip} 0& 0& 0& -(23/9)& 0& -(20/3)& 0& 0
 \\\noalign{\medskip} 0& 0& -(17/18)& 0& -(23/9)& 0& -(20/3)& 0
\\\noalign{\medskip}0& -(1/108)& 0& -(55/54)& 0& -(14/9)& 0& 0
\\\noalign{\medskip}0& 0& -(1/108)& 0& -(55/54)& 0& -(14/9)& 0
\\\noalign{\medskip}0& 1/648& 0& 1/324& 0& -(20/27)& 0& 0
\\\noalign{\medskip}0& 0& 1/648& 0& 1/324& 0& -(20/27)& 0
\\\noalign{\medskip}0& -8& 0& -64& 0& -96& 0& 0
\\\noalign{\medskip}0& 0& 0& 0& 0& 0& 0& 0
\\\noalign{\medskip}0& 0& 0& 0& 0& 0& 0& 0
\\\noalign{\medskip}-2& 0& -8& 0& -32& 0& -32& 0 \end{array} \right)
$$ }

\vspace{2pc}
$$  \frac{\partial T}{\partial x_1}(0,-1/2, 1/2)=$$
{\tiny $$  \left (\begin {array}{ccccccccccccc}0& -(7/4)& 0& -5& 0& -7& 0& 0
 \\\noalign{\medskip}0& 0& 0& 0& -(155/36)& 0& -(35/6)& 0
 \\\noalign{\medskip}0& 1/72& 0& -(19/12)& 0& -(10/3)& 0& 0
 \\\noalign{\medskip}-(1/432)& 0& -(1/72)& 0& -(367/216)& 0& -(127/36)& 0
 \\\noalign{\medskip}0& -(1/432)& 0& -(1/72)& 0& -(10/9)& 0& 0
 \\\noalign{\medskip}1/2592& 0& 1/432& 0& 7/1296& 0& -(233/
  216)& 0
 \\\noalign{\medskip}0& 1/2592& 0& 1/432& 0& 5/27& 0& 0
 \\\noalign{\medskip}-2& -4& -12& -24& 0& 0& 0& 0
 \\\noalign{\medskip}0& 0& 8& 0& 32& 0& 32& 0
 \\\noalign{\medskip}0& 4& 0& 8& 0& 0& 0& 0
 \\\noalign{\medskip}0& 0& 0& 0& 0& 0& 0& 0
\end {array} \right)
$$ }
$$  \frac{\partial T}{\partial x_2}(0,-1/2, 1/2)=$$
{\tiny $$ \left (\begin{array}{ccccccccccccc}
\\\noalign{\medskip}-1& 0& -(3/2)& 0& -3& 0& 0& 0
\\\noalign{\medskip}0& 0& 0& -(4/3)& 0& -3& 0& 0
\\\noalign{\medskip}0& 0& -(35/36)& 0& -(4/3)& 0& -3& 0
\\\noalign{\medskip}0& -(1/216)& 0& -1& 0& -(5/6)& 0& 0
\\\noalign{\medskip}0& 0& -(1/216)& 0& -1& 0& -(5/
  6)& 0
\\\noalign{\medskip}0& 1/1296& 0& 0& 0& -(31/
  36)& 0& 0
\\\noalign{\medskip}0& 0& 1/1296& 0& 0& 0& -(31/
  36)& 0
\\\noalign{\medskip}0& -8& -12& -64& -100& -96& -168& 0
\\\noalign{\medskip}0& 4& 0& 8& 0& 0& \
0& 0
\\\noalign{\medskip}2& 0& 0& 0& 0& 0& 0& 0
\\\noalign{\medskip}0& 0& 0& 0& 0& 0& 0& 0\end {array} \right)
$$  }

\newpage

{\epsfysize=35mm \epsfbox{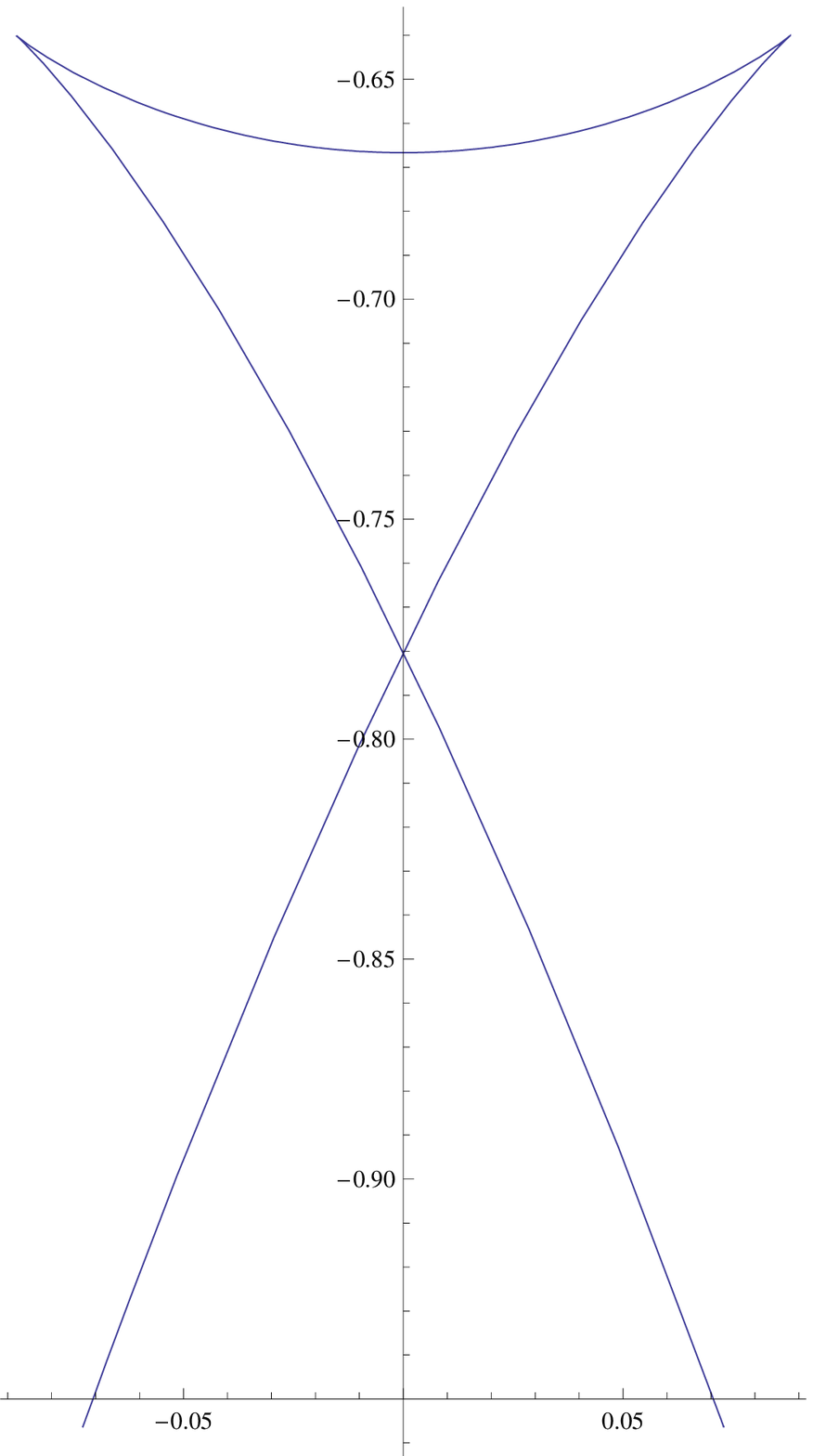} t=2/3}, \;\;\;\; {\epsfysize=35mm \epsfbox{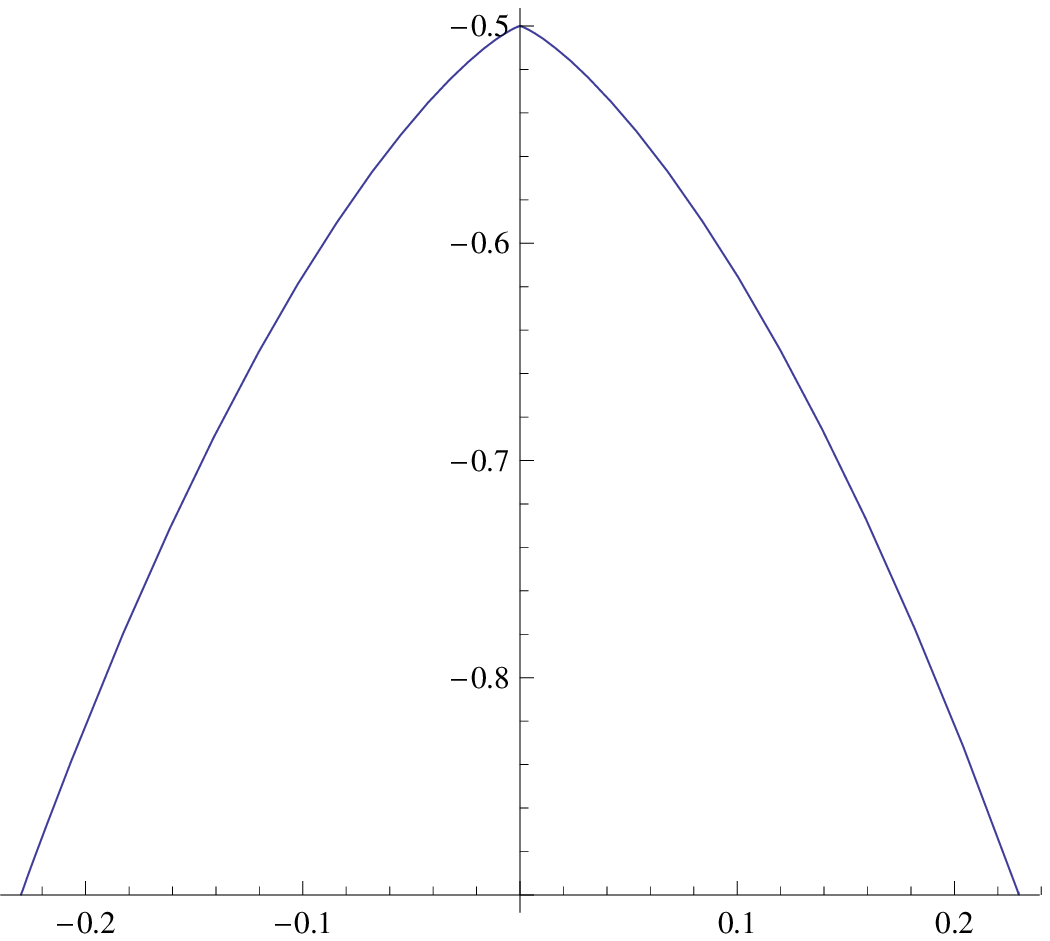} t=1/2, \;\;{\bf real \;wave\; fronts} }

\vspace{1pc}

\;\;\; {\epsfysize=35mm \epsfbox{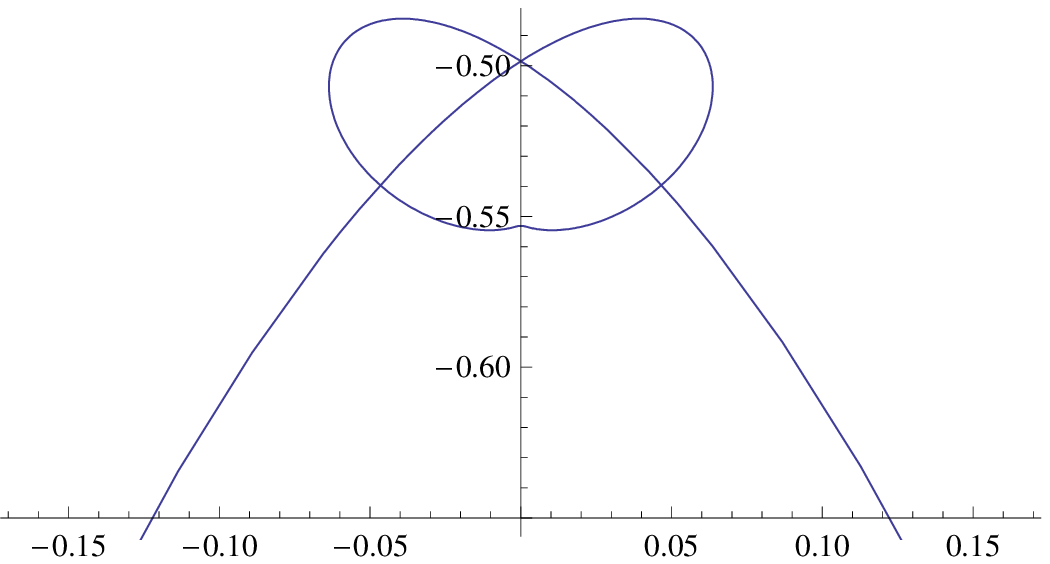} t=0.55, \;\; $\Im z$= 1/4 }

\;\;\; {\epsfysize=35mm \epsfbox{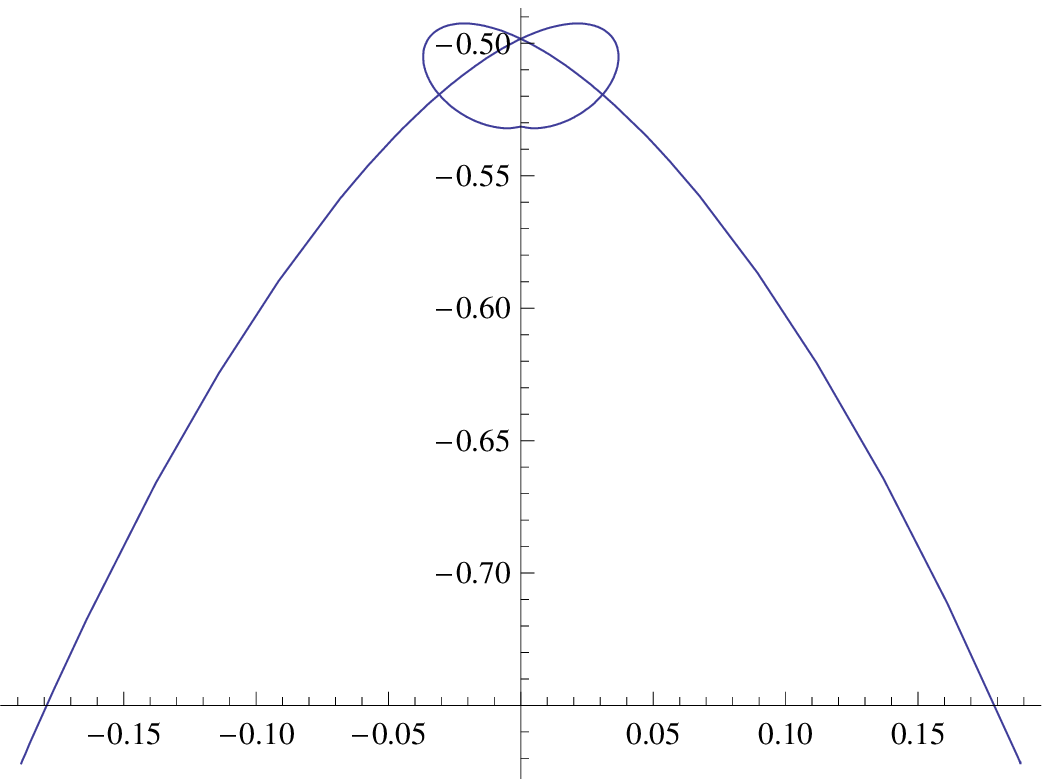} t=0.53,  \;\; $\Im z$= 0.175 }   

\;\;\;  {\epsfysize=35mm \epsfbox{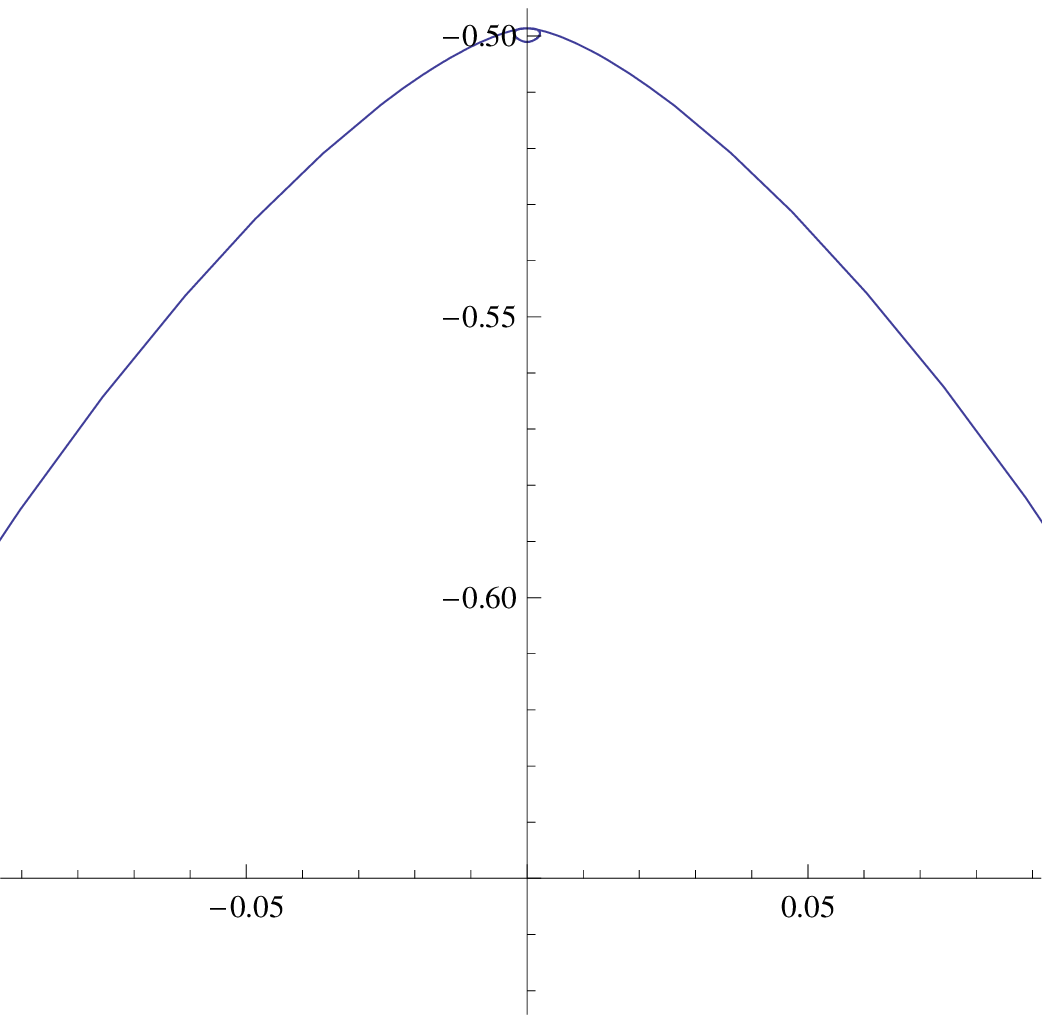} t=0.501, \;\;  $\Im z$ = 0.1 }

\vspace{2pc}
{\bf Figures} showing successive degeneration of wave fronts. 
The first two figures show real wave fronts at the moments $t=2/3, 1/2.$

The following three figures show real sections of the complex wave fronts at the moments $t=0.55, 0.53, 0.501$
that started from the initial front with corresponding to fixed $\Im z$ values. One can verify that in the neighbourhood
of a $A_3$ focal point no trefoil shaped figure appears. The wave front remains to be more or less a parabola shape figure.

\newpage

Now we shall show that $8\times 11$ matrix of $(5.6)$ attains the maximal rank $8$ except at $(\bC^3,0)\ni(\xi_1,\xi_2,\tau)=(0,0,0)$. For this purpose, we select $8 \times 8$ minors without a common factor vanishing at $(\xi_1,\xi_2,\tau)=(0,0,0).$
Here $[i,j,k,...]$ stands for a $8 \times 8$ minor corresponding to the $i,j,k,...$-th rows of $(5.6).$

{\small $[2, 3, 4, 5, 6, 9, 10, 11]=-(1/629856)(1 + 2 \tau) (-26375 \xi_1^4 + 2736000 \xi_1^6 + 
    ...- 7587840 \xi_2 \tau^5)$}

{\small $[1, 2, 3, 4, 7, 9, 10, 11]=-(1/1259712)\xi_1 (1 + 2 \tau) (76350 \xi_1^2 - 
    164788000 \xi_1^4 + ... + 
    1221525504 \tau^5)$}
   
{\small $[2, 3, 4, 7, 8, 9, 10, 11]=
    (1/104976)(1 + 2 \tau) (-22905 \xi_1^2 + 
    ... + 32102400 \tau^6)$}

According to a calculation by SINGULAR ~\cite{GrPf}, system of algebraic equations $ [ 2,$ $3,$ $4,$ $5, 6, 9, 10, 11]$ $=[1, 2, 3, 4, 7, 9, 10, 11]$ $ = [2, 3,$ $4,$ $7,$ $8, 9, 10, 11] = 0$ has
$(\xi_1,\xi_2,\tau)=(0,0,0)$ as an isolated solution with multiplicity 12. This means that the sufficient condition of Proposition~\ref{proposition44} is satisfied.}
This means that  $(0,-1/2, 1/2)$ is an isolated point on $BW$ $\subset \iota^{-1}(D_\varphi)$ with the property $(5.6).$
Upshot is the almost freeness of the big wave front germ at the focal point after Proposition ~\ref{proposition42}.

In summary we gave two proofs to the statement.
\begin{proposition}
The germ of the big wave front $BW$ at the focal point $(x_1, x_2,t)=(0,-1/2a, 1/2a)$ 
defines a free divisor if $ a \not= 1.$
If $a=1$ it defines an almost free divisor germ at the focal point $(x_1, x_2,t)=(0,-1/2, 1/2)$. 
\end{proposition}

{\bf 2. Wave propagation in the 3 dimensional space}

 Now we consider the following initial wave front in the $3-$dimensional space,
$Y:=\{(z,u) \in \bC^{2}: -\frac{1}{2}(k_1 z_1^2 +k_2 z_2^2)+u=0  \},$ i.e. $F (z)= -\frac{1}{2}(k_1 z_1^2 +k_2 z_2^2)$ for $0<k_1 < k_2.$ 
In this case our phase function has the following expression
$$\Psi(x,t,z)=  (-x_3 + k_1 x_1 z_1 + k_2 x_2 z_2 - 1/2 (k_1 z_1^2 + k_2 z_2^2))^2 - 
 t^2 (1 + k_1^2 z_1^2 + k_2^2 z_2^2),$$
$$ =-t^2 + x_3^2 - k_1^2 x_1 z_1^3 + (k_1^2 z_1^4)/4 - 2 k_2 x_3 (x_2 - z_2) z_2 $$
$$-  k_2^2 t^2 z_2^2 - k_2 x_3 z_2^2 + k_2^2 (x_2 - z_2)^2 z_2^2 + k_2^2 (x_2 - z_2) z_2^3 + (k_2^2 z_2^4)/4 $$
$$+ z_1^2 (-k_1^2 t^2 + k_1^2 x_1^2 + k_1 x_3 - k_1 k_2 (x_2 - z_2) z_2 - 
    1/2 k_1 k_2 z_2^2) $$
$$+ z_1 (-2 k_1 x_1 x_3 + 2 k_1 k_2 x_1 (x_2 - z_2) z_2 + k_1 k_2 x_1 z_2^2) \eqno(5.7)$$ 
It is easy to see that the point $(x_1,x_2,x_3,t) =(0,0,1/k_1, 1/k_1)$ is a focal point with
a singular point $(z,u) =(0,0)$ and the Milnor number 
$\mu(0)=3.$ We have the following tame polynomial,
$$\Psi(0,0,1/k_1, 1/k_1,z)=(
k_1^4 z_1^4 + 4 k_1 k_2 z_2^2 - 4 k_2^2 z_2^2 + 2 k_1^3 k_2 z_1^2 z_2^2 + 
 k_1^2 k_2^2 z_2^4)/4 k_1^2.$$
 As a matter of fact, the polynomial $\Psi(0,0,1/k_1, 1/k_1,z)$ satisfies the criterion on the presence of $A_3$
singularity at the origin mentioned in  ~\cite{Hase}, Theorem  2.2, (2). The situation is the same at another  focal point $(x_1,x_2,x_3,t) =(0,0,1/k_2, 1/k_2).$ 
 The quotient ring $(1.5)$ for this $\Psi(0,0,1/k_1, 1/k_1,z)$ has dimension $\mu=5$.
 
 We can choose  
$$ \{e_1,e_2,e_3,e_4,e_5\} = \{  1,z_1,z_1^2, z_2,z_2^2 \}$$
as the basis $(2.8).$  In view of $(5.7),$ we introduce additional deformation monomials
$e_6= z_1*z_2$,$e_7=z_2^3$,
$e_8=z_1^3$, $e_9=z_1^2*z_2$, $e_{10}= z_1*z_2^2$ together with the entries of the mapping 
$\iota,$
$$s_1=-t^2 + x_3^2, s_2=-2 k_1 x_1 x_3,s_3= -k_1^2 t^2+k_1^2 x_1^2 + k_1 x_3, s_4=-2 k_2 x_2 x_3$$
$$ s_5= - k_2^2 t^2 + k_2^2 x_2^2 + k_2 x_3, s_6=2 k_1 k_2 x_1 x_2$$  
$$s_7=-k_2^2 x_2,s_8= -k_1^2 x_1, s_{9}=-k_1 k_2 x_2,s_{10}=-k_1 k_2 x_1. $$

By direct calculation with the aid of  SINGULAR ~\cite{GrPf}, we can verify 
 $$  dim_{\bC}\frac{\bC[z]}{d_z(\Psi(0,0,1/k_1, 1/k_1,z) + \sum_{i=1}^6 s_i e_i)\bC[z]} =5,$$
 while 
 $$ dim_{\bC}\frac{\bC[z]}{d_z(\Psi(0,0,1/k_1, 1/k_1,z) + \sum_{i=1}^6 s_i e_i + s_j e_j)\bC[z]} =7,$$
for $j=7,8,9,10.$  
If $s \in \bC^{10}$ satisfies the condition
 $$ s_7*z2^3+s_8*z_1^3+s_9*z_1^2*z_2+s_{10}*z_1*z_2^2 =(\alpha z_1 + \beta z_2)(k_1 z_1^2 +k_2 z_2^2), \eqno(5.8)$$
 for some $(\alpha, \beta) \in \bC^2, $ $\varphi(z,s)$ has the same singularity at the infinity ($A_1 + A_1$) for every such $s.$ In fact our $\Psi(x,t,z)$ is exactly of this form with $ (\alpha, \beta) =(-k_1 x_1,- k_2 x_2).$ This means that for
 $s \in \bC^{10}$ under the condition $(5.8)$, the global total Milnor number $\mu$ is equal to $5$.

 Additionally we remark here that for $\Psi(0,0,1/k_1, 1/k_1,z) + \sum_{i=1}^6 s_i e_i + s_j e_j$,$j=7,8,9,10$ the jump ($=2$) of Milnor number infinity takes place as $s_j \rightarrow 0.$  This illustrates the upper semi-continuity of the Milnor number at infinity.
 
In summary, we can choose a constructible set 
$$U = \{s' \in \bC^{9}; k_1 s_7- k_2 s_9=0, k_1 s_{10} - k_2 s_8=0 \}$$ 
with dimension $7$ for which Lemma~\ref{lemma2} applies.
 This implies that the Assumption I,(i) is satisfied with $\nu=8$.

The above discussion proves that the image of the mapping $ \iota(\bC^4) \subset \bC^{10}$ 
 is contained in a constructible set $\bC \times U$  where the value of the matrix $\Sigma(s)$ is well-defined at each point 
$s \in \bC \times U$.
Therefore
$$ dim_{\bC}\frac{\bC[z]}{d_z(\Psi(x,t,z))\bC[z]} =5, $$
 for every $(x,t) \in \bC^{4}.$ This means that the Assumption I,(ii) is satisfied.  
 
 At the focal point $(x_1,x_2,x_3,t) =(0,0,1/k_1, 1/k_1)$ the matrix $\iota^\ast(\Sigma)$
has the following form with rank $3$
$$
\left (\begin {array}{ccccccccccccccc}0& 0& 0& 0& -k_2 (k_1 - k_2)/2 k_1^2& 0 &0 &0
\\\noalign{\medskip}0& 0& 0& 0& 0& 0& 0& 0
 \\\noalign{\medskip}0& 0& 0& 0& 0& 0& 0& 0
 \\\noalign{\medskip}0& 0& 0& (k_1 - k_2)^2/k_1^4& 0& 0& 0& 0
 \\\noalign{\medskip}0& 0& 0& 0& (k_1 - k_2)^2/k_1^4& 0& 0& 0
 \\\noalign{\medskip}0& 0& 0& 0& 0& -1& 0& 0
 \\\noalign{\medskip}0& 0& 0& 0& 0&  0&-1& 0
 \\\noalign{\medskip}0& 0& 0& 0& 0&  0& 0&-1
 \end {array} \right)
$$

Together with the data    
     $$d_{x,t} \iota(0,0,1/k_1, 1/k_1)= $$
 {\small   $$ \left (\begin {array}{ccccccccccccc}0& -2& 0& 0& 0& 0& 0& -k_1^2& 0& -k_1 k_2&
     \\\noalign{\medskip}0& 0& 0& -2 k_2/k_1& 0& 0& -k_2^2& 0& -k_1 k_2& 0
     \\\noalign{\medskip}2/k_1& 0& k_1& 0& k_2& 0& 0& 0& 0& 0
     \\\noalign{\medskip}-2/k_1& 0& -2 k_1& 0& -2 k_2^2/k_1& 0& 0& 0& 0& 0\end {array} \right)$$}
 we see that the $rank\; T(0,0,1/k_1, 1/k_1)= 8 \geq \nu =8.$
By virtue of the Proposition ~\ref{proposition42}, we have the following.

\begin{proposition} The wave front contained in the discriminantal loci of $(5.7)$
defines a free divisor germ in the neighbourhood of
the focal point $(0,0,1/k_1, 1/k_1).$
\end{proposition}

\noindent

\begin{flushleft}
 \begin{minipage}[t]{8cm}
  \begin{center}
{\footnotesize   
Susumu Tanab\'e\\
 Previous address:  Kumamoto University\\
Kurokami 2-39-1,\\
 Kumamoto,860-8555,
Japan\\
\vspace{1pc} 
Current address:Galatasaray University\\
\c{C}{\i}ra$\rm\breve{g}$an cad. 36,\\
Be\c{s}ikta\c{s}, \.{I}stanbul, 34357,
Turkey\\
   Emails:tanabesusumu@hotmail.com, tanabe@gsu.edu.tr}
\end{center}
\end{minipage}\hfill
\end{flushleft}

\end{document}